\documentclass[11pt,a4paper]{article}
\usepackage{amsfonts,amssymb}

\usepackage{amscd}

\vspace{0.5cm}

 4

\font\erm=cmr8

\usepackage{graphicx}

\author{M.~Dziemia\'nczuk}
\title{Counting Bipartite, k-Colored Multi  
and Directed Acyclic Multi Graphs Through F-nomial coefficients}

\newtheorem{defn}{Definition}

\newtheorem{theoremn}{Theorem}
\newtheorem{observen}{Observation}
\newtheorem{notation}{Notation}
\newtheorem{fact}{Fact}

\newtheorem{corollary}{Corollary}

\newtheorem{lemma}{Lemma}

\newcommand{\layer}[2]{\langle\Phi_{#1} \! \to\! \Phi_{#2}\rangle}

\newcommand{\fnomial}[2]{ {{#1} \choose {#2}}_{\!\!F} }
\newcommand{\fnomialF}[3]{ {{#1} \choose {#2}}_{\!\!#3} }

\newcommand{\sm}{\scriptstyle}

\begin{document}

\begin{center}
\noindent {\Large \textsc{Counting Bipartite, k-Colored \\
and Directed Acyclic Multi Graphs Through F-nomial coefficients}}  \\ 

\vspace{0.5cm}

\noindent {Maciej Dziemia\'nczuk}

\vspace{0.5cm}

\noindent {\erm Gda\'nsk University Student, the Institute of Computer Science}

\noindent {\erm PL-80-952 Gda\'nsk, st. Wita Stwosza 57, Poland}

\noindent {\erm e-mail: maciek.ciupa@gmail.com}
\end{center}

\vspace{1cm}

\noindent \textbf{Abstract} 

\noindent F-nomial coefficients encompass  among others well-known binomial coefficients or Gaussian coefficients that count subsets of finite set and subspaces of finite vector space respectively. Here, the so called F-cobweb tiling sequences $N(\alpha)$ are considered. For such specific sequences a new interpretation with respect to Kwa\'sniewski general combinatorial interpretation of $F$-nomial coefficients is unearhed.
 
Namely, for tiling sequences $F = N(\alpha)$ the $F$-nomial coefficients are equal to the number of labeled special bipartite multigraphs denoted here as $\alpha$-multigraphs $G(\alpha,n,k)$. 

An explicit relation between the number of $k$-colored $\alpha$-multigraphs and multi $N(\alpha)$ -nomial coefficients is established. 
We also prove that the unsigned values of the first row of inversion matrix for $N(\alpha)$ -nomial coefficients considered here are equal to the numbers of directed acyclic $\alpha$-multigraphs with $n$ nodes.

\vspace{0.4cm}
\noindent AMS Classification Numbers: 05A19 , 11B39, 15A09.


\vspace{0.2cm}
\noindent Affiliated to The Internet Gian-Carlo Polish Seminar:

\noindent \emph{http://ii.uwb.edu.pl/akk/sem/sem\_rota.htm}

\section{Introduction}

The notation from \cite{akk1,akk2} is being here taken for granted.

\vspace{0.2cm}
\noindent \textbf{Comment 1}

\noindent  For the mnemonic efficiency of Kwa\'sniewski  up-side-down notation see Appendix in \cite{akk4}  and references therein and consult recent \cite{akk5}, \cite{akk7}. 
\noindent With this Kwa\'sniewski "upside down notation" inspired by Gauss  and applying reasonings almost just repeated with "$k_F$" numbers   replacing  $k$ - natural numbers one gets in the spirit of Knuth \cite{knuth}  clean results also in this report.
\noindent And more ad "upside down notation":  concerning Gauss and Knuth - see remarks in \cite{knuth} also on Gaussian  binomial coefficients.

\begin{defn}[\cite{akk1,akk2}]
Let any $F$-cobweb admissible sequence then $F$-nomial coefficients are defined as follows
$$
	\fnomial{n}{k} = \frac{n_F!}{k_F!(n-k)_F!} 
	= \frac{n_F\cdot(n-1)_F\cdot ...\cdot(n-k+1)_F}{1_F\cdot 2_F\cdot ... \cdot k_F}
	= \frac{n^{\underline{k}}_F}{k_F!}
$$
\noindent while $n,k\in \mathbb{N}$ and $0_F! = n^{\underline{0}}_F = 1$.
\end{defn}

\begin{notation}
Let us denote by $\mathcal{T}_\lambda$ a family of natural numbers' valued sequences $F\equiv\{n_F\}_{n\geq 0}$ constituted by $n$-th coefficients of the generating function $F(x)$ expansion i.e. $n_F = [x^n]F(x)$  (in Wilf's notation \cite{wilf}), where
\begin{equation}
	F(x) = 1_F\cdot \frac{x}{(1 - \alpha x)(1 - \beta x)}
\end{equation}
\noindent for $1_F\in\mathbb{N}$ and $\alpha,\beta\in \mathbb{R}$.
\end{notation}

It was shown in  \cite{md3} that any $F\in\mathcal{T}_\lambda$ is \emph{cobweb-tiling} sequence.

\begin{notation}
$N(\alpha) = F = \{n_F\}_{n\geq 0}\in\mathcal{T}_\lambda$ means that $\alpha=\beta$ and $1_F=1$.
\end{notation}

\noindent $N(\alpha)$ is a cobweb tiling sequence, of course.

\vspace{0.4cm}
Let us recall some of the sequence $F \equiv \{n_F\}_{n\geq 0} = N(\alpha)$ properties following \cite{md3}. 
Let at first $n>1$ and $1_F\in\mathbb{N}$ be given, then 

\begin{equation}\label{eqAlpha}
	n_F = 1_F \cdot n \cdot \alpha^{n-1}
\end{equation}

\noindent Then the following recurrence relation for any $m,k\in \mathbb{N}$ takes place

$$
	n_F = (k+m)_F = \alpha^m k_F + \alpha^k m_F
$$

\noindent Hence the corresponding $F$-nomial coefficients do satisfy

$$
	\fnomial{n}{k} = \alpha^m \fnomial{n-1}{k-1} + \alpha^k \fnomial{n-1}{k}
$$

\noindent while $\fnomial{n}{n} = \fnomial{n}{0} = 1$.

\section {Combinatorial interpretation of $N(\alpha)$ tiling sequences $N(\alpha)$-nomials}

\subsection{$F$-cobweb admissible sequence}

At first, let us refer to the base joint Kwa\'sniewski combinatorial interpretation of $F$-nomial coefficients of all at once cobweb-admissible sequences. Then we recall some special properties of the sequence $N(\alpha)$ following \cite{md3,md4}.

\begin{fact}\label{fact:1}
For $F$-cobweb admissible sequences $F$-nomial coefficient $\fnomial{n}{k}$ is the number of max-disjoint \emph{equipotent} copies $\sigma P_{n-k}$ of the layer $\layer{k+1}{n}$. (\cite{akk1,akk2} and references therein)
\end{fact}

\textbf{Specifically important}: $F$-cobweb tiling sequences' $F$-nomial coefficients are the numbers of max-disjoint equipotent copies $\sigma P_{n-k}$ - \textbf{the layer $\layer{k+1}{n}$ is tiled with}.

\vspace{0.4cm}
Because of that we consider only just these sequences $N(\alpha)\in \mathcal{T}_\lambda$, in what follows.

\vspace{0.4cm}
\noindent \textbf{Recall.} 
The sequences $F=N(\alpha)\in \mathcal{T}_\lambda$  have several combinatorial interpretations. For example, if $a=2$ then $n_F$ equals the number of edges of an $n$-dimensional hypercube and the number of 132-avoiding permutations of $[n+2]$ containing exactly one $123$ pattern (sequence $A001787$ \cite{oeis})

\vspace{0.4cm}
\noindent \textbf{Examples of the sequences $N(\alpha) = \{ n_{N(\alpha)}\}_{n\geq{0}}$}

\begin{enumerate}
\item $N(1) = 0, 1, 2, 3, 4, 5, 6, 7, 8, 9, 10, ... $
\item $N(2) = 0, 1, 4, 12, 32, 80, 192, 448, 1024, 2304, 5120, ...  $
\item $N(3) = 0, 1, 6, 27, 108, 405, 1458, 5103, 17496, 59049, 196830, ...  $
\item $N(4) = 0, 1, 8, 48, 256, 1280, 6144, 28672, 131072, 589824, 2621440,...   $

\end{enumerate}

\noindent \textbf{Observe.} $F= N(\alpha)$. Let $\alpha \in \mathbb{N}$ be given, then $n_F=n\cdot \alpha^{n-1}$ i.e. it is equal to the number of sequences $(s_1,s_2,...,s_n)$ such that one of terms is equal to zero i.e. $\exists! k : (1\leq k \leq n \wedge s_k=0)$, and the rest of them $s_j\in [\alpha]$.

\subsection{Counting bipartite and $K$-colored $\alpha$-multigraphs}

Recall. Labeled bipartite $\alpha$-multigraph $G(\alpha,n,k)$ is a bipartite graph with $n$ vertices ($k$ of them is in one two disjoint vertices' sets) with multiedges, such that any two vertices might be connected by at most $(\alpha-1)$ edges. We define $K$-colored $\alpha$-multigraph in a similar way.

\vspace{0.4cm}
\noindent Note. If $\alpha=1$ then the graph $G(\alpha,n,k) = (V_1 \cup V_2,E)$ has no edges i.e. $E=\emptyset$. Therefore it might be considered as a $k$-subset of $n$-set i.e. $G(1,n,k) \equiv V_1 \subseteq [n]\wedge |V_1|=k$ . Then labeled $K$-colored $1$-multigraph is considered as a partition of set $[n]$ into $k$-nonempty blocks where each of them represents a set of vertices' indices of $G$ with the same color.

\begin{observen}\label{obs:1}
The number of labeled bipartite $\alpha$-multigraphs $G(\alpha,n,k)$ denoted by $\beta_{\alpha,n,k}$ is 
\begin{equation}
	\beta_{\alpha,n,k} = {n \choose k}\cdot \alpha^{k(n-k)}
\end{equation}
\end{observen}

\noindent {\it{\textbf{Proof.}}}
Take any $\alpha\in\mathbb{N}$. If $\alpha=1$ then we have no edges between nodes, therefore we need to count all $k$-subsets of vertices set $n$. If $\alpha>1$ then there is additionally $\alpha^{k(n-k)}$  possibilities to create at most $(\alpha-1)$ edges between any two vertices from disjoint sets $V_1$ and $V_2$ $\blacksquare$

\vspace{0.4cm}
\noindent Note, $\beta_{2,n,k}$ is equal to the number of bipartite graphs $G(2,n,k)$ of $n$ nodes where $k$ of them belong to one of disjoint vertices' sets.

\begin{corollary}\label{cor:1}
Let $F$ be a cobweb tiling sequence $N(\alpha) \in \mathcal{T}_{\lambda}$, such that $1_F=1$. Then the $F$-nomial coefficient is equal to the number $\beta_{a,n,k}$ of labeled bipartite graphs $G(\alpha,n,k)$
\begin{equation}\label{eqFib}
	\fnomial{n}{k} = \beta_{\alpha,n,k}
\end{equation}
\noindent for $n,k\in\mathbb{N}\cup\{0\}$.
\end{corollary}

\noindent {\it{\textbf{Proof.}}}
Consider the sequence $F\equiv\{n_F\}_{n\geq 0} = N(\alpha)$ i.e. such that $1_F=1$ and $n_F = n\cdot \alpha^{n-1}$.
We only need to show that the number of bipartite $\alpha$-multigraphs $G(\alpha,n,k)$ is $\beta_{\alpha,n,k}$ (Observation \ref{obs:1}). For that to see take any $n,k\in\mathbb{N}$, and just check that 
$$
\fnomial{n}{k} = \frac{n\alpha^{n-1} \cdot (n-1)\alpha^{n-2} \cdot ... \cdot (n-k+1)\alpha^{n-k}}{1\cdot\alpha^0 \cdot 2\alpha^1 \cdot ... \cdot k\alpha^{k-1}} = 
$$
$$
	= \frac{n^{\underline{k}}}{k!} \cdot \frac{\alpha^{k\frac{2n-k-1}{2}}}{\alpha^{k\frac{k-1}{2}}}
	= {n \choose k} \cdot \alpha^{k(n-k)} \blacksquare
$$

\vspace{0.4cm}

In another words, $\fnomial{n}{k} = \fnomialF{n}{k}{N(\alpha)}$ is equal to the number of $2$-colored $\alpha$-multigraphs with $n$ vertices, while $k$ of them are colored by one color and $n-k$ by another one.

\vspace{0.4cm}
\noindent \textbf{Note.} Consider $F=N(2)$. The number $\gamma_{n,2}$ (see: \cite{finch}) of all 2-colored graphs is equal to

\begin{equation}
	\gamma_{n,2} = \sum_{k\geq 0} {n\choose k} \cdot 2^{k(n-k)} = \sum_{k\geq 0} \fnomialF{n}{k}{N(2)}
\end{equation}

In general, for $\alpha\in\mathbb{N}$ the sum $\sum_{n\geq 0}\fnomialF{n}{k}{N(\alpha)}$ is equal to the number of all $2$-colored $\alpha$-multigraphs.

\vspace{0.4cm}
\noindent \textbf{Comment 2} \emph{On F-binomiality and 2-colored $\alpha$-multigraphs.}

\noindent (Due to  A. Krzysztof Kwa\'sniewski, see \cite{akk6} and consult for notations also references [7,8,9] therein).

In view of the final Remark in \cite{akk6},  the combinatorics fundamental logarithmic Fib-Binomial Formula (\cite{akk6}, Section 4)

$$
	\phi_n^{(t)}(x +_F a) \equiv \left[ \mathrm{exp}\{a\partial_F\}\phi_n^{(t)} \right](x)
	= \sum_{k\geq 0}\left[\!\! \begin{array}{c} n \\ k \end{array}\!\!\right]_F \phi_{n-k}^{(t)}(a)x^k
$$
$$
	t=0,1; \ \ \ |x| < a; \ \ \ n\in \mathbb{Z}
$$

\noindent may be considered as $F$-Binomial for any natural numbers valued sequence $F$ with $F_0=1$  (the class considered in \cite{akk6} is much broader). For special $F$-sequences known as $F$-cobweb posets admissible sequences the $F$-nomial coefficients   

\begin{equation}\label{eq:hybrid}
	\left[\!\! \begin{array}{c} n \\ k \end{array}\!\!\right]_F = \fnomial{n}{k}
\end{equation}

\noindent for $n,k\geq 0$ acquire Kwa\'sniewski joint combinatorial interpretation (Fact \ref{fact:1}).

\vspace{0.2cm}
Now put in $F$-binomial formula above $t=0$ and $a=x=1$ and pay attention to - that according to the definition of $F$-hybrid binomial coefficients one has (\ref{eq:hybrid}) for $n,k\geq 0$. Then we get (known at least since Morgan Ward famous Calculus of Sequences \cite{ward}) the clean - appealing in Kwa\'sniewski notation formula:

\begin{equation}
	\left( 1 +_F 1 \right)^n \equiv \sum_{n\geq 0}\fnomial{n}{k}
\end{equation}

In particular and rephrasing: the number $\gamma_{n,2}$ of all 2-colored graphs as in \cite{finch} is now equal to

$$
	\gamma_{n,2} = \sum_{k\geq 0} {n\choose k} \cdot 2^{k(n-k)} = \sum_{k\geq 0} \fnomialF{n}{k}{N(2)} \equiv \left( 1 +_F 1 \right)^n
$$

\vspace{0.2cm}

Let us show up some values of $\{\gamma_{n,2}\}_{i\geq 0}$ (sequence $A047863$ \cite{oeis})

\vspace{0.2cm}
\noindent $\{\gamma_{n,2}\}_{i\geq 0} =  {1, 2, 6, 26, 162, 1442, 18306, 330626, 8488962, ... }$

\vspace{0.4cm}
\noindent And here is the matrix $\mathbf{M}=[m_{ij}]$, $m_{ij} = \fnomialF{i}{j}{N(2)}$ of the sequence $N(2)$, where $0\leq i,j\leq 7$.

$$
	M = \left[
	\begin{array}{llllllll}
	1 \\
	1 &  1 \\
	1 &  4   & 1 \\
	1 &  12  & 12    & 1 \\
	1 &  32  & 96    & 32     & 1 \\
	1 &  80  & 640   & 640    & 80     & 1  \\
	1 &  192 & 3840  & 10240  & 3840   & 192   & 1 \\
	1 &  448 & 21504 & 143360 & 143360 & 21504 & 448  &  1
	\end{array}
	\right]
$$

Due to Corollary \ref{cor:1} and thanks to general properties of $F$-cobweb admissible tiling sequences we infer families of identities. 
For example recurrence relation for the number $\beta_{\alpha,n,k}$ of labeled bipartite $\alpha$-multigraphs $G(\alpha,n,k)$ is given "for free" and it reads

\begin{equation}
	\beta_{\alpha,n,k} = \fnomialF{n}{k}{N(\alpha)}
	= \alpha^{n-k}\cdot \fnomialF{n-1}{k-1}{N(\alpha)} + \alpha^k\cdot \fnomialF{n-1}{k}{N(\alpha)}
\end{equation}

It has got combinatorial proof in cobweb posets language as it is the case with all cobweb-tiling sequences \cite{md3}.

It might be expressed also in bipartite graphs terms with the use of a standard counting rule. We fix the last $n$-th vertex and separate family of all graphs $G(\alpha,n,k)$ into two disjoint classes depending on that, where the vertex is assigned.

\begin{observen}
Let $c_\alpha(\vec{b})$ be the number of labeled $k$-colored $\alpha$-multigraphs $G(\alpha,\vec{b})$ with $n$ vertices, where $\vec{b} = \langle b_1,b_2,...,b_k \rangle$ such that $b_1$ vertices are colored by first color, next $b_2$ vertices by another one and so on. Then $c_\alpha(\vec{b})$ is 

\begin{equation}
	c_\alpha(\vec{b}) = {n \choose {b_1,b_2,...,b_k}} \cdot \alpha^{\frac{1}{2}\left(n^2 - b_1^2 - b_2^2 - ... - b_k^2\right)}
\end{equation}

\noindent while $n = b_1 + b_2 + ... + b_k$.
\end{observen}

\noindent {\it{\textbf{Proof.}}}
Take any such vector $\vec{b} = \langle b_1,b_2,...,b_k \rangle$. The coloring of $n$ vertices might be chosen on ${n \choose {b_1,b_2,...,b_k}}$  ways. Additionally, if $\alpha = 1$ then a graph $G(\alpha,\vec{b})$ has no edges but if $\alpha > 1$, then 
there might be created at most $(\alpha - 1)$ edges on $\alpha^{b_i\cdot b_j}$ ways between any two vertices from disjoint vertices' sets $V_i$, $V_j$ where $i\neq j$. Therefore the overall number of all possibilities is
$$
	\alpha^{\sum_{1\leq i < j \leq k} b_i\cdot b_j}
	= \alpha^{b_1 b_2 + b_1 b_3 + ... + b_{k-1} b_k}
	= \alpha^{\frac{1}{2}\left(n^2 - b_1^2 - b_2^2 - ... - b_k^2\right)}
$$

\noindent Hence the thesis $\blacksquare$

\begin{corollary}\label{cor:2}
Let $F$ be a cobweb tiling sequence $N(\alpha) \in \mathcal{T}_\lambda$ such that $1_F=1$. 
Then the multi $F$-nomial coefficient is equal to the number $c_\alpha(\vec{b})$ of labeled bipartite $\alpha$-multi graphs $G(\alpha,\vec{b})$ i.e.

\begin{equation}
	\fnomial{n}{b_1,b_2,...,b_k} = c_\alpha(\vec{b})
\end{equation}

\noindent where $\vec{b} = \langle b_1,b_2,...,b_k \rangle$ and $b_1 + b_2 + ... + b_k = n$
\end{corollary}

\noindent \textbf{Note.}
Let $\gamma_{\alpha,n,k}$ be the number of all $k$-colored $\alpha$-multigraphs with $n$ vertices. Then 

\begin{equation}
	\gamma_{\alpha,n,k} = \sum_{{b_1+...+b_k=n \atop b_1,...,b_k\geq 0}} \fnomialF{n}{b_1,b_2,...,b_k}{N(\alpha)}
\end{equation}

The case of $\alpha=2$ i.e. when $G(\alpha,\vec{b})$ is a $k$-colored graph without multiple edges was already considered in \cite{finch}.

\section {Counting labeled directed acyclic $\alpha$-multigraphs ($\alpha$-DAGs) }

\begin{defn}
A directed acyclic graph with $\alpha$-multiple edges i.e. any two vertices might be connected by at most $(\alpha-1)$ directed edges is called acyclic $\alpha$-multi digraph ($\alpha$-DAG for short).
\end{defn}

\begin{lemma}\label{lem:1}
Let $A_\alpha(n)$ denotes the number of acyclic $\alpha$-multi digraphs ($\alpha$-DAGs) with $n$ labeled nodes. Then for $n\in\mathbf{N}$

\begin{equation}\label{eq:dags}
	A_\alpha(n) = \sum_{k\geq 1} (-1)^{k+1} \fnomialF{n}{k}{N(\alpha)}\cdot A_\alpha(n-k)
\end{equation}

\noindent while $A_\alpha(0) = 1$ and $\alpha\in\mathbb{N}$.
\end{lemma}

\noindent {\it{\textbf{Proof.}}}
The main idea of the proof comes from \cite{robin} (see also \cite{robin0,robin1}) where particular case of $\alpha=2$ is considered with the help of inclusion-exclusion principle. One shows that any directed acyclic multi-graph with no cyclic paths has at least one vertex with in-degree equal to zero (such vertices are so-called \emph{out-points} \cite{robin}).

\vspace{0.2cm}
Take $\alpha\in\mathbb{N}$ and a graph $\alpha$-DAG with $n\in\mathbb{N}$ nodes.
Denote by $X_i$ a family of $\alpha$-DAGs, such that $i$-th point is an out-point for $1\leq i \leq n$.
Therefore $A_\alpha(n) = \left| \bigcup_{i=1}^{n} X_i \right|$  and from inclusion-exclusion principle 

$$
	A_\alpha(n) = \left| \bigcup_{i=1}^{n} X_i \right| 
	= \sum_{k=1}^{n} (-1)^{k+1} \sum_{1\leq b_1<...<b_k\leq n} \left| X_{b_1} \cap ... \cap X_{b_k} \right|
$$

Let us consider the number $k\in[n]$ of out-points. We can label them on ${n\choose k}$ ways. Next, there are $\alpha^{k(n-k)}$ possibilities to eventually create $\alpha$-multiple edges from these $k$ points to the rest $(n-k)$ ones for which we can create $A_\alpha(n-k)$ $\alpha$-DAGs, thus 
$$
	 \sum_{1\leq b_1<...<b_k\leq n} \left| X_{b_1} \cap ... \cap X_{b_k} \right| 
	= {n \choose k}\alpha^{k(n-k)}\cdot A_\alpha(n-k)
$$

\noindent Therefore

$$
	A_\alpha(n) = \sum_{k=1}^{n} (-1)^{k+1} {n \choose k} \alpha^{k(n-k)} \cdot A_\alpha(n-k)
$$

\noindent and according to Corollary \ref{cor:1}

$$
	A_\alpha(n) = \sum_{k=1}^{n} (-1)^{k+1} \fnomialF{n}{k}{N(\alpha)} \cdot A_\alpha(n-k)
$$

\noindent Hence the thesis $\blacksquare$

\begin{fact}
Let $F$-cobweb admissible sequence be given. Then an inversion formula for $F$-nomial coefficients derived in \cite{md4} is of the form 

$$
	\fnomial{n}{k}^{-1} = \fnomial{n}{k} \fnomial{n-k}{0}^{-1}, \ \ \ \ \fnomial{n}{n}^{-1} = 1
$$
\begin{equation}
	\fnomial{n}{0}^{-1} = \sum_{s=1}^k (-1)^s \sum_{k_1+...+k_s=n \atop k_1,...,k_s\geq 1} \fnomial{n}{k_1,k_2,...,k_s}
\end{equation}

\end{fact}

Here is the  inversion matrix $\mathbf{M}^{-1} \equiv \fnomial{i}{j}^{-1}$ of matrix $\mathbf{M}$  from previous section example i.e. for $F=N(2)$.

$$
	M = \left[
	\begin{array}{llllllll}
		\sm 1                                                               \\
		\sm -1          & \sm  1                                                    \\
		\sm 3           & \sm -4         & \sm 1                                         \\
		\sm -25         & \sm  36        & \sm -12        & \sm 1                              \\
		\sm 543         & \sm -800       & \sm 288        & \sm -32      & \sm 1                     \\
		\sm -29281      & \sm  43440     & \sm -16000     & \sm 1920     & \sm -80      & \sm 1            \\
		\sm 3781503     & \sm -5621952   & \sm 2085120    & \sm -256000  & \sm 11520    & \sm -192  & \sm 1      \\
	\end{array}
	\right]
$$

\begin{theoremn}
Let $F$ be a cobweb tiling sequence $F=N(\alpha)\in\mathcal{T}_\lambda$, such that $1_F=1$ and let $A_\alpha(n)$ denotes the number of labeled acyclic $\alpha$-multi digraphs with $n$ vertices. Then for $n\in\mathbb{N}$

\begin{equation}
	A_\alpha(n) = (-1)^n \fnomial{n}{0}^{-1} = \left| \fnomial{n}{0}^{-1} \right|
\end{equation}

\noindent where $A_\alpha(0) = 1$ and $\fnomial{n}{k}^{-1}$ stays for an inversion matrix of $F$-nomial coefficients.
\end{theoremn}

\noindent {\it{\textbf{Proof.}}}
Take any $n\in\mathbb{N}$. If Lemma \ref{lem:1} is taken into account, then
$$
	A_\alpha(n) = \sum_{1\leq k_1\leq n} (-1)^{k_1+1} \fnomial{n}{k_1}\cdot A_\alpha(n-k_1)
$$

\noindent While expanding the above we set new sums' variables as $k_1,k_2,...,k_n$ and there are at the most $n$ variables $k_s$ with each of them equal one, according to the conditions $k_1+...+k_n = n$, with $k_i \geq 0$

\begin{eqnarray}
\lefteqn{
	A_\alpha(n) = 
	\sum_{k_1=1}^{n} (-1)^{k_1+1}\fnomial{n}{k_1} 
	\sum_{k_2=1}^{n-k_1} (-1)^{k_2+1}\fnomial{n-k_1}{k_2} 
	 \cdots
} & &
\nonumber\\
& & \cdots  \sum_{k_n=1}^{n-k_1-...-k_{n-1}} (-1)^{k_n+1}\fnomial{n-k_1-...-k_{n-1}}{k_n} 
	\cdot A_\alpha(0) 
\end{eqnarray}

\noindent and consequently

$$
	A_\alpha(n) = \sum_{I}
	(-1)^{k_1+k_2+...+k_n + S}
	\fnomial{n}{k_1,k_2,...,k_n}
$$

\noindent where $S$ is the number of variables $k_1,k_2,...,k_n$ with positive value and
$$
	I =  \left\{ \begin{array}{l}
	\sm 0\leq k_1 \leq n \\
	\sm 0\leq k_2 \leq n-k_1 \\
	\sm ... \\
	\sm 0\leq k_n \leq n-k_1-...-k_{n-1}
	\end{array} \right.
$$

\noindent Now, let us rearrange the sum into two summations  as follows 

$$
	A_\alpha(n) = \sum_{s=1}^n  (-1)^s  \sum_{k_1+...+k_s=n \atop k_1,...,k_s\geq 1} \!\!\!\!\!(-1)^n \fnomial{n}{k_1,...,k_s} 
	= (-1)^n \fnomial{n}{0}^{-1}
$$

\noindent The value of $A_\alpha(n)$ is positive for any natural $n$, hence the thesis $\blacksquare$

\vspace{0.8cm}
\noindent \textbf{Acknowledgements}

\vspace{0.2cm}
I would like to thank Professor A. Krzysztof Kwa\'sniewski - who initiated
my interest in his cobweb poset concept - for his very helpful comments and
improvements of this note.


\end{document}